\documentclass[12pt,leqno]{article}
\usepackage[latin1]{inputenc}
\usepackage{amsmath,amsfonts,amssymb,amsthm,amscd}
\date{}

\newcommand\Z{\mathbb{Z}}

\newcommand\A{\mathbb{A}}

\newcommand\C{\mathbb{C}}

\newcommand\R{\mathbb{R}}

\newcommand\Q{\mathbb{Q}}

\newcommand\fa{\mathfrak{a}}

\newcommand\fb{\mathfrak{b}}
\newcommand\fd{\mathfrak{d}}

\newcommand\pp{\mathfrak{p}}

\newcommand\fq{\mathfrak{q}}

\newcommand\ff{\mathfrak{f}}

\newcommand\Hh{\mathcal{H}}

\newcommand\Cc{\mathcal{C}}
\newcommand\Mm{\mathcal{M}}

\newcommand\CPp{\mathcal{P}}

\newcommand\Ll{\mathcal{L}}

\newcommand\Ff{\mathcal{F}}

\newcommand\Frob{\mathrm{Frob}}
\newcommand{\ba}{\backslash}
\newcommand\rg{\rightarrow}
\newcommand\lgr{\longrightarrow}   
\newcommand\Gal{\mathrm{Gal}}

\newcommand\Sup{\mathrm{Sup}}

\newcommand{\ga}{\mathbf{A}}

\newcommand{\gf}{\mathbf{f}}

\newcommand{\gD}{\mathbf{D}}

\newcommand{\Inf}{\mathrm{Inf}}

\newcommand{\nequiv}{\equiv\kern-4.5mm/}

\def\adots{\mathinner{\mkern1mu\raise1pt\vbox{\kern7pt\hbox{.}}
\mkern2mu\raise4pt\hbox{.}
\mkern2mu\raise7pt\hbox{.}\mkern1mu}}

\newcommand{\intt}{\int\kern-4.5mm-}

\newtheorem{monlem}{{Lemma}}[section]

\newtheorem{maprop}{{Proposition}}[section]
\newtheorem{montheo}{{Theorem}}[section]

\numberwithin{equation}{section}

\begin{document}

\title{New applications of the Mellin transform\\
to automorphic $L$--functions}

\author{Laurent Clozel}

\maketitle

\noindent \textbf{Introduction}

In an earlier paper \cite{Cl}, and its Appendix, written with Peter Sarnak, we have obtained universal lower bounds on certain quadratic integrals of automorphic $L$--functions. For instance, if $\pi$ is a cuspidal unitary representation of $GL(m,\A_\Q)$, and $L(s)=L(s,\pi)$:
\begin{equation}
\int_{-\infty}
^{+\infty} \Big|\frac{L(1/2+it)}{1/2+it}\Big| ^2 dt > \frac{\pi}{2}.\end{equation}
Cf \cite[Theorem D]{Cl}. In this article, we follow some questions that arose naturally in this context.

The first one was suggested by a remark of Sarnak, according to which we cannot have such universal bounds for short intervals, e.g. for the integral on $[-1,1]$ in (0.1); at least such cannot be obtained if $m$ is allowed to vary. See the Introduction to \cite{Cl}, \S~2.

The argument does not succeed if $m$ is fixed; neither did  the author succeed in finding, for $m$ fixed, an absolute bound on a short interval. This remains an interesting problem.

In Chapter 2, we obtain a universal lower bound, for $m$ fixed, for the integral on an interval $[-A\log C, A\log C]$ where $C$ is the analytic conductor. Even then, we could not obtain a fixed lower bound. Rather, we prove that this integral is larger than $c(\log C)^{-1/2}$ where $c>0$ is an absolute constant. This may not be optimal, but it is commensurable with the Lindel\"of conjecture. (We hope that an analytic number theorist familiar with the use of the Mellin transform, as in the proof of the approximate functional equation, will be able substantially to improve this result.) Moreover, we can shift the ordinate and obtain a general result on an interval $[X-T,X+T]$ where $X$ is arbitrary and $T$ is of the order of $\log X$. See Theorem~2.2, Theorem~2.3. The proof relies on a theorem of Molteni~\cite{Mo}, further improved by Xiannan Li \cite{Li}.

In Chapter 3, we follow a lead from \cite[\S 3.3]{Cl}. There we considered Vinogradov's Conjecture on the order of the first quadratic non--residue and we showed that it followed directly, via the Mellin transform, from the Lindel\"{o}f conjecture (including in the $q$--aspect) for the associated Dirichlet $L$--function. (This may have been well--known to experts).

A similar problem is, given a non trivial representation $\rho$ of a Galois group $\Gal(E/F)$ of number fields, to determine ``the'' first prime $\pp$ of $F$ --- i.e. one of smallest norm --- where $\rho$ is unramified and $\rho(\Frob_\pp)\not=1$, $\Frob_\pp$ being a Frobenius  element. The proof of \cite{Cl} extends naturally. We have treated first the case where $\rho$ is a one--dimensional character, in Chapter~3, \S~1--3. See Theorem~3.1. When $\rho$ is non--Abelian, the proof is more delicate and leads us to introduce the unramified variant $\Ll(s,\rho)$ of the Artin $L$--function $L(s,\rho)$. See Theorem~3.2.

In \S 3.6, we show that our estimate is, in some cases, better that the known unconditional estimates. (The relevant result here is the recent one of Zaman \cite{Z2}) However we have to assume that $F$ is a Galois extension of $\Q$ and, more crucially, that $L(s,\rho)$ is holomorphic. Note that we do not have to assume that $\rho$ is associated to an automorphic representation.

In  Chapter 4, we have reviewed some odds and ends concerning the results and arguments of \cite{Cl}. In particular we point out that they apply to Rankin $L$--functions --- not only to the standard $L$--functions of \cite{Cl}; and we develop the remarks made in  \cite{Cl} about the relation between estimates on the summation function $A_0(x)$ --- see \S~2.1--- and subconvexity for~$L(s,\pi)$.

Finally, in Chapter 1, we have recalled some ``well--known'' results concerning the growth of $L$--functions in the critical strip, in term of the analytic conductor. Fortunately we could rely on the very clear exposition of Harcos \cite{Ha}. For the convenience of the reader, we also collect some formulas for special functions.

We use notation standard in analytic number theory, in particular Landau's symbols $\ll,O(\ )$, indexed when we want to specify the dependence of the implicit constants. We use $f \prec g$ ($f,g$ positive functions) for $f /g \rightarrow 0$. We use $A$ for an absolute constant (in the context), not always the same in different occurrences.

\vskip2mm

\textit{Acknowledgement.-} I thank Farrell Brumley,  Jesse Thorner and Asif Zaman for providing useful references, and Peter Sarnak for suggesting that my lower bounds should be compared with the well-known bound of Ramachandra. I also thank the Fondation Simone et Cino del Duca for financial support.

\section{Majorations and formulas}

\subsection{} In this preliminary chapter we have grouped together the majorations and formulas which will be used in the text, and crucially in using an ``absolute'' version, due to Molteni, of the Friedlander--Iwanie\v c estimate. We briefly recall our set--up, which is that of \cite{Cl}. We  consider a unitary cuspidal representation $\pi$ of $GL(m,\A)$, or a product $\pi=\pi_1\times \pi_2\times \cdots\times \pi_r$ (parabolic induction) of cuspidal representations $\pi_i$ of $GL(m_i,\A)$, $m=\sum m_i$. Then
$$
L(s,\pi) = \prod_{i=1}^r L(s,\pi_i).
$$
is the standard $L$--function. The completed $L$--function
$$
\Lambda(s,\pi) = D^{s/2}\ L(s,\pi_\infty)\ L(s,\pi),
$$
where $D$ is the conductor, an integer $\ge 1$, satisfies a functional equation
$$
\Lambda(s,\pi) = \varepsilon(\pi) \Lambda(1-s,\tilde{\pi}).
$$
For details and a review of the other $L$--functions to which these theorems apply, see \cite[\S~2.1]{Cl}.

We do not assume, as we did in \cite{Cl}, that $\pi_\infty$ is self--dual. The functional equation can than be written
\begin{equation}
L(1-s,\pi) = \varepsilon(\pi)\gamma(s)L(s,\tilde{\pi})
\end{equation}
with
$$
\gamma(s) = (\pi^{-m}D)^{s-1/2} \frac{c(\tilde{\pi}_{\infty})\Gamma(s,\tilde{\pi}_\infty)}{c(\pi_{\infty})\Gamma(1-s,\pi_\infty)}
$$
and 
$$
\Gamma(s,\pi_\infty) = \prod_{j=1}^m \Gamma\Big(\frac{s+c_j}{2}\Big)
$$
with
$$
Re(c_j) \ge \frac{1}{m^2+1}-\frac{1}{2}\qquad \cite{LRS}.
$$

Following Iwanie\v c and Sarnak, we associate to $\pi$ its conductor $D=D(\pi)$ and its analytic conductor
\begin{equation}
C=C(\pi) = D(\pi) \prod_{j=1}^m (2+|c_j|)
\end{equation}
as well as
$$
C(s) = C(\pi,s) = C(\pi)\ (1+|s|)^m.
$$
(There is a finer version of $C(s)$ ; cf \cite[p.~95]{IK} for a discussion of this.) 
We will need uniform estimates, in terms of these data, for $\gamma(s)$ and $L(s,\pi)$ in a strip $Re(s)\in ]-\varepsilon,1+\varepsilon[$.

\subsection{} We recall the known bounds on $\gamma(s)$ and $L(s)$. For $\gamma(s)$ a uniform bound is derived by Harcos \cite{Ha}. See however the corrections in \cite{Ha2}, in particular \cite[(3)]{Ha2} which corrects \cite[3.22]{Ha}\footnote{I thank Farrell Brumley for this reference.}.

\vspace{2mm}

{\monlem  (Harcos) {For} $\sigma >\frac{1}{m}-\frac{1}{m^2+1}$, $Re(s)=\sigma$,}
$$
\gamma(s) \overset{}{\underset{\sigma}{\ll}} C(s)^{\sigma-1/2}.
$$

We now assume the Ramanujan conjecture for $\pi$. Fix $\varepsilon>0$ (small.) Since $\pi_f$ is tempered, we have $L(s,\tilde{\pi}) \ll 1$ for $Re(s)=1+\varepsilon$, uniformly for $m$ fixed. For such a value of $s$, $\gamma(s)\ll C(s)^{1/2+\varepsilon}$. The functional equation then implies
$$
L(s,\pi) \ll C(s)^{1/2+\varepsilon}
$$
for $Re(s)=-\varepsilon$.

We deduce from the Phragm\'en--Lindel\"of principle:

\vspace{2mm}

{\maprop Assume $\pi_f$ tempered. For $-\varepsilon\le Re(s) \le 1+\varepsilon$,}
$$
L(s,\pi) \ll_\varepsilon C(s)^{\frac{1-\sigma}{2}+\varepsilon}.
$$
Cf. Iwanie\v{c}--Kowalski \cite[p. 100]{IK}.

\subsection{} Here we insert a few remarks on the analytic conductor. It is not true that, for $A>0$ fixed, there exist a finite number of representations $\pi$ such that $C(\pi)\le A$. Indeed, if $m=1$, and $\pi=\chi$ is a Dirichlet character of conductor $D$, twisted by $|\,|^a$, $(a\in i\R)$, the analytic conductor is $D(2+|a|)$. This phenomenon persists in higher rank $m$, but is essentially due to the center $GL(1)$ of $GL(m)$.

\vspace{2mm}

{\maprop For any $A>0$, there exist a finite number of cuspidal representations $\pi$ of $GL(m,\A)$ whose central character $\omega$ verifies $\omega|_{\R_+^\times}=1$ and such that $C(\pi)<A$.}

\vspace{2mm}

Since $D(\pi) < C(\pi)$ there is a finite number of possibilities for $D(\pi)$. Let $S$ denote the connected component of 1 in the center $Z(\R)\cong \R^\times$ of $GL(m,\R)$. Consider the space $\Ll=L_{cusp}^2(S\, GL(m,\Q)\ba GL(m,\A))$ of $L^2$--cusp forms invariant by $S$. The Hecke algebra $\bigotimes\limits_v \Hh_v$ (where $\Hh_\infty = C_c^\infty(GL(n,\R))$ and, for $v=p$ finite, $\Hh_p$ is the algebra of compactly supported smooth functions) acts on $L_{cusp}^2$ by right translations, and it is well--known that this action is trace--class. Fix $D\ge 1$, and let $K(D) \subset GL(m,\A_f)$ be the  congruence subgroup defined by Jacquet, Piaterskii-Shapiro and Shalika \cite[\S 5, Th\'{e}or\`{e}me]{JPSS}, and $\varphi_D\in \bigotimes\limits_p \Hh_p$ its characteristic function. Then $\varphi_D$ projects $\Ll$ onto its subspace composed of the $K(D)$--invariants in the cuspidal representations $\pi$ (with $\omega_\pi |_{\R_+^\times}=1)$. This is an infinite sum of representations $\rho$ of $GL(m,\R)$; a function $\varphi_\infty \in C_c^\infty(GL(m,\R))$ acts on it by a trace--class operator. In particular, any compact subset of the unitary dual of $GL(m,\R)$ contains only a finite number of such representations. But the set of representations $\rho$ such that $C_\infty(\rho)= \prod (2+|c_j|)\le A$ is compact. (See \cite[\S~2.2]{Cl} for the relation between the $c_j$ and the Langlands parameters of~$\pi_\infty$.)

\vspace{2mm}

Assume again $\pi$ cuspidal. For $a\in i \,\R$, let $\pi[a]=\pi\otimes |\det|^a$. We can choose $a$  such that $\pi[a]$ has trivial central character on $\R_+^\times$. In view of this, it is useful to know the relation between $C(\pi)$ and $C(\pi[a])$. 

Again we refer to \cite[\S~2.2]{Cl}. The representation $\pi_\infty$ is associated to a sum of real characters
$$
\nu(x) = (sgn x)^\varepsilon|x|^c \qquad (x\in \R_+^\times);\varepsilon=0,1)
$$
and of characters of $\C^\times$:
$$
\mu(z) =z^p(\bar{z})^q,\ p-q\in \Z.
$$
(The result of Luo--Rudnick--Sarnak \cite{LRS} implies simple bounds on $Re(c)$,\break $Re(p+q)$.)

The datum $c_j$ associated to $\nu$ is $c+\varepsilon$; the data $c_j$ associated to $\mu$ are
$$
\begin{array}{lll}
(q,q+1) &\mathrm{if} &Re(p-q) <0\\
(p,p+1) &\mathrm{if} &Re(p-q)>0.
\end{array}
$$

If we twist by $|\ |^a$, the $c_j$ are transformed into $c_j+a$. Therefore
\begin{equation}
C(\pi[a]) \le C(\pi)(1+|a|)^m.
\end{equation}

We also remark that if $C(\pi)=2^m$, $D(\pi)=1$ and the $c_j$ are equal to 0. Under this assumption, the number of $\pi$ is therefore finite.

Finally, it is of interest to know when $\omega|_{\R_+^\times} =1$ in terms of the $c_j$. However the relation is not direct. The real characters $\nu_j$ contribute $c_j=c+\varepsilon$, where $\nu_j(x)=x^c$, $x>0$; $A$ complex character $\mu(z)$ determines a representation $\pi_2$ of $GL(2,\R)$ of central character $x^{p+q}$ $(x>0)$, and $c_j$, $c_{j'}$ equal to $(q,q+1)$ or $(p,p+1)$. One checks that $c_j+c_{j'}=p+q+n$ where $n=|p-q|>0$ is the ramification of $\pi_2$. Therefore the central character of $\pi_\infty$, on $\R_+^\times$, is $x^A$,
$$
A=\sum c_j -Ram(\pi_\infty),
$$
where
$$
Ram(\pi_\infty) = \sum \varepsilon+\sum n
$$
is the ramification of $\pi_\infty$. To obtain a representation with $\omega|_{\R_+^\times}=1$, we have to twist by $|\ det|^{-A/m}$. The new conductor is (very roughly) bounded by 
$$
C(\pi) (1+|A/m|)^m.
$$ 

\subsection{} For the reader's convenience we recall a few facts on special functions. We will use the functions
$$
\begin{array}{rl}
Erf(x) &=\displaystyle\int_0^x e^{-t^2}dt,\\
\noalign{\vskip2mm}
Erfc(x) &=\displaystyle\int_x^\infty e^{-t^2} dt = \frac{1}{2}\sqrt{\pi} - Erf(x).
\end{array}
$$

See \cite[p. 147]{Bate}. In particular we have the  asymptotic expansion of $Erfc(x)$ for $x\gg0$:
$$
Erfc(x) =\frac{1}{2}e^{-x^2}(x^{-1}+O(x^{-3})).
$$
We will also need an estimate for
$$
I(a,T) = \int_T^\infty e^{-t^2} t^a dt\qquad (T>0).
$$
This is best computed as an incomplete gamma function. With
$$
\begin{array}{rl}
\Gamma
(a,x) &= \displaystyle \int_x^\infty e^{-t} t^a \frac{dt}{t}\,,\\
\noalign{\vskip2mm}
I(a,T) &=\displaystyle \int_T^\infty e^{-t^2} t^{a+1} \frac{dt}{t}\\
\noalign{\vskip2mm}
&=\displaystyle\frac{1}{2} \int_{T^2}^\infty e^{-s}s^{\frac{a+1}{2}} \frac{ds}{s}\\
\noalign{\vskip2mm}
&=\displaystyle\frac{1}{2}\Gamma\Big(\frac{a+1}{2},T^2\Big)
\end{array}
$$
and the asymptotic expansion \cite[p. 135]{Bate}
$$
\Gamma(a,x) = x^{a-1}e^{-x}(1+O(\frac{1}{x}))
$$
now yields
$$
I(a,T) = \frac{1}{2}T^{a-1} e^{-T^2}(1+O(T^{-2})).
$$

\section{The quadratic integral of $\frac{L(1/2+it)}{1/2+it}$ on short intervals}

\subsection{} In this section we assume $m$ fixed; we consider $\pi$ verifying the assumptions in (1.1). We write 
$$
\begin{array}{rl}
L(s)&=L(s,\pi),\\

\noalign{\vskip2mm}
L(s,\pi) &=\displaystyle \sum a_nn^{-s},\\
\noalign{\vskip2mm}
A_0(x) &= \displaystyle\sum_{n\le x} a_n\quad (x\ge 1).
\end{array}
$$
We assume given $\xi$, $\nu>0$ and consider only $\pi$ such that $\nu\ge \sigma_c(\pi)$ where $\sigma_c(\pi)$ is the abscissa of convergence of $L(s,\pi)$. We assume that we have, for $x\ge 1$ and any $\varepsilon$, an estimate 
\begin{equation}
A_0(x) \ll_\varepsilon C^\xi \, x^{\nu+\varepsilon}
\end{equation}
where the implicit constant depends only on $\varepsilon$, and $C=C(\pi)$. For a suitable function $f$ on $\R_+^\times$, we define
$$
\widetilde {\Mm}f(s) = \int_0^\infty f(x) x^{-s}\frac{dx}{x}.
$$

{\monlem (cf. $\cite[\S~3.3]{Cl})$ For $Re(s)>\nu$, the integral defining $\widetilde{\Mm} A_0(s)$ is absolutely convergent and $\widetilde{\Mm} A_0(s) = \frac{L(s)}{s}$.}

\vspace{2mm}

Indeed
$$
\begin{array}{l}
\hbox to 1cm{}\displaystyle\int_0^X \Big(\sum_{n\le x} a_n\Big) n^{-s}\frac{dx}{x}\\
\noalign{\vskip2mm}
=\displaystyle\sum_{n\le X} a_n \int_n^X  x^{-s-1}dx\\
\noalign{\vskip2mm}
= \displaystyle\frac{1}{s}\sum_{n\le X}a_n n^{-s} - \frac{1}{s}\Big(\sum_{n\le X}a_n\Big) X^{-s}.
\end{array}
$$
Since $Re(s) >\sigma_c(\pi)$, the first sum converges to $L(s,\pi)$ for $X\rg \infty$. The second is dominated by $X^{\nu+\varepsilon/2} X^{-\nu-\varepsilon} \rg 0$.

For $\sigma>\nu$, we then have
\begin{equation}
\int_1^\infty x^{-\sigma} A_0(x) x^{-it} \frac{dx}{x} = \frac{L(\sigma+it)}{\sigma+it}.
\end{equation}

In order to obtain a minoration  of the $L^2$--norm of $\frac{L(\frac{1}{2}+it)}{1/2+it}$ on an interval $[-T,T]$, we consider its scalar product with a Gaussian. Thus let, for $\alpha>0$:
$$
g_\alpha(t)= e^{-\pi\alpha t^2}.
$$
We extend $g_\alpha$ to a function of $s$, equal to $g_\alpha(t)$ for $s=1/2+it$. Thus
$$
G_\alpha(s)=e^{\pi\alpha(s-1/2)^2}.
$$

This is a function of rapid decrease in $t$ $(s=\sigma+it)$, uniformly in any vertical strip. In  particular, for~$\sigma>\nu$:
\begin{equation}
\int_\sigma \frac{L(s)}{s} G_\alpha(s)ds=i \int_{-\infty}^{+\infty} \frac{L(1/2+it)}{1/2+it} g_\alpha(t)dt
\end{equation}
since $L(s)$ and $G_\alpha(s)$ have no poles. We estimate the left--hand side using the Fourier transform. We have by (2.2)
$$
\int_0^\infty e^{-\sigma X} A_0 (e^X) e^{-itX}dX = \frac{L(\sigma+it)}{\sigma+it},
$$
i.e. 
$$
\frac{L(\sigma+it)}{\sigma+it}= \Ff(e^{-\sigma X}A_0(e^X))
$$
where
$$
\Ff h(t) = \int_{-\infty}^{+\infty} h(X) e^{-itX}dX.
$$

The inverse transformation is
\begin{equation}
\Ff^{-1}k(X) = \frac{1}{2\pi} \int_{-\infty}^{+\infty} k(t) e^{itX}dt.
\end{equation}{
The scalar product formula is
$$
\int k_1(X) \overline{k_2(X)} dX =\frac{1}{2\pi} \int h_1(t) \bar{h}_2(t)dt.
$$
However, (2.3) is a bilinear product, and then:
$$
\int h_1(t) h_2(t)dt = 2\pi \int k_1(X) k_2(-X)dX.
$$
Therefore the left--hand side of (2.3) is equal to the product of $i=\sqrt{-1}$ and of
\begin{equation}
2\pi \int_0^\infty e^{-\sigma X} A_0(e^X) \hat{G}_\alpha(-X)dX
\end{equation}
where $\hat{G}_{\alpha}(X)$ is given by (2.4) applied to $k(t)=G_\alpha(\sigma+it)$. Now
\begin{equation}
G_\alpha(\sigma+it) = e^{\pi\alpha(\sigma-1/2)^2} e^{-\pi\alpha t^2} e^{2\pi\alpha(\sigma-1/2)it}.
\end{equation}
We have 
$$
\Ff^{-1} k(X) = \frac{1}{2\pi}\Ff_0^{-1}k(\frac{X}{2\pi})
$$ 
where $\Ff_0$, $\Ff_0^{-1}$ denote the usual Fourier transforms, normalised by $2\pi$. Now
$$
\begin{array}{rl}
\Ff_0^{-1}
(f(t)e^{2i\pi\beta t}) &=\Ff_0^{-1} f(X+\beta),\\
\Ff_0^{-1}(e^{-\pi\alpha t^2}) &=\dfrac{1}{\sqrt{\alpha}} e^{-\frac{\pi}{\alpha}X^2}
\end{array}
$$
so $\Ff_0^{-1}G_\alpha)(X)$ is, by (2.6), equal to
$$
\begin{array}{c}
e^{\pi\alpha(\sigma-1/2)^2} \dfrac{1}{\sqrt{\alpha}}e^{-\frac{\pi}{\alpha}(X+\alpha(\sigma-1/2))^2}\\
= \dfrac{1}{\sqrt{\alpha}} e^{-\frac{\pi}{\alpha}X^2} e^{-2\pi X(\sigma-1/2)},
\end{array}
$$
and
$$
\Ff^{-1}G_\alpha(-X) = \frac{1}{2\pi\sqrt{\alpha}} e^{-\frac{X^2}{4\pi\alpha}} e^{(\sigma-1/2)X}.
$$

\subsection{} We now consider the integral (2.5), equal to
\begin{equation}
\begin{array}{c}
\displaystyle \frac{1}{\sqrt{\alpha}} \int_0^\infty e^{-\sigma X} A_0(e^X) e^{-\frac{X^2}{4\pi \alpha}} e^{(\sigma-1/2)X} dX\\
= \displaystyle \frac{1}{\sqrt{\alpha}} \int_0^\infty A_0(e^X) e^{-\frac{X^2}{4\pi \alpha}} e^{-\frac{X}{2}} dX.
\end{array}
\end{equation}
(It is, as it should be in view of the translation of complex integrals, independent of $\sigma$.)

For $1\le x\le 2$, $A_0(x)=1$. We first obtain a lower bound for
$$
I_1 := \frac{1}{\sqrt{\alpha}} \int_0^{\log 2} e^{-\frac{X^2}{4\pi \alpha}} e^{-\frac{X}{2}}dX.
$$
In the sequel, Landau's symbols $O(\ ),\ll,\ldots,$ unindexed, are used when the implicit constants are absolute. We have
$$
I_1 \geq \sqrt{1/2} \frac{1}{\sqrt{\alpha}}
 \int_0^{\log2} e^{-\frac{X^2}{4\pi\alpha}} dX.
 $$
 We set $X=\sqrt{4\pi\alpha}Y$; the integral is then
$$
\begin{array}{rl}
&2\sqrt{\pi}\displaystyle \int_0^{\log 2/2\sqrt{\pi\alpha}} e^{-Y^2}dY\\
= &\displaystyle 2\sqrt{\pi} Erf \Big(\frac{\log 2}{2\sqrt{\pi\alpha}}\Big)\\
= &2\sqrt{\pi} \displaystyle \Big(\frac{1}{2} \sqrt{\pi} -Erfc\Big(\frac{\log 2}{2\sqrt{\pi\alpha}}\Big)\Big).
\end{array}
$$
For $x\rg \infty$, $Erfc(x)=\frac{1}{2x} e^{-x^2}(1+O(\frac{1}{x^2}))$. For small $\alpha$, then,
$$
I_1 \ge \pi \sqrt{1/2} -O\Big(\sqrt{\alpha}\ e^{-\frac{\log^2 2}{4\pi \alpha}}\Big)
$$
and $I_1 \ge \frac{\pi}{2}:=2c$ for $\alpha$ sufficiently small.

We now have to estimate the remainder, dominated by
$$
I_2 = \frac{1}{\sqrt{\alpha}} \int_{\log 2}^\infty C^\xi \ e^{(\nu+\varepsilon-1/2)X} e^{-\frac{X^2}{4\pi\alpha}}dX.
$$
Write $\theta=\nu+\varepsilon-1/2$, and set
$$
Y=X-2\pi \theta\alpha.
$$
The exponential term in the integrand is then
$$
\exp\Big(-\frac{Y^2}{4\pi\alpha}\Big) e^{\pi\theta^2\alpha}.
$$
We can neglect the constant since $\alpha$ will be small. Thus $I_2$ is dominated by
$$
\begin{array}{rl}
I_3 &=\displaystyle \frac{1}{\sqrt{\alpha}} \int_\ell^\infty C^\xi e^{-\frac{1}{4\pi\alpha}X^2}dX,\\
&\ell=\log 2-2\pi\theta\alpha.
\end{array}
$$

for small $\alpha$, so the part relative to $[\ell,\log 2]$ is dominated by $C^{\xi}\sqrt{\alpha}\cdot e^{-a_1/\alpha}$ with $a_1>0$ fixed. We now consider
$$
I_4 = C^\xi \int_{\log 2}^\infty \frac{1}{\sqrt{\alpha}} e^{-\frac{1}{4\pi\alpha}X^2} dX.
$$

The integral is equal to $2\sqrt{\pi}$ $Erfc(\frac{\log 2}{2\sqrt{\pi\alpha}})$ (see the formulas in \S~1); it admits an expression
\begin{equation}
\sqrt{\pi} \ e^{-\frac{\log^2 2}{4\pi\alpha}} \Big(\frac{2\sqrt{\pi\alpha}}{\log 2} +O(\alpha^{3/2})\Big)
\end{equation}

for small $\alpha$. This (multiplied by $C^\xi$) is of the same order as the integral on the small segment.

We want to ensure that this is dominated by $I_1$, which is implied by
$$
\sqrt{\alpha} \ e^{-\frac{a_2}{\alpha}} C^\xi \le a_3
$$
where $a_3$ is a small constant, i.e.
$$
\frac{1}{2}\log \alpha - \frac{a_2}{\alpha} + \xi\log C \le -a_4
$$
$(a_4>0)$. We may assume $\alpha\le 1$, so we seek 
\begin{equation}
\frac{a_2}{\alpha} \ge a_4 + \xi \log C.
\end{equation}
Since $C\ge 2^m$, $\xi\log C>a_5>0$ and
$$
\frac{1}{a_4+\xi \log C} > \frac{1}{N\xi \log C}
$$
for $N$ such that $a_4<(N-1)\xi\log C$, so $N$ is determined by constants independent of~$\pi$. 

Thus (2.9) is verified if
\begin{equation}
\alpha
\le \frac{a_5}{\log C}.
\end{equation}

We summarise the result

\vspace{2mm}

{\maprop Assume $\pi$ cuspidal, and $\sigma_c(\pi)\le \nu$. There exist positive constants $c$, $b_1$ (depending only on $\xi$, $\nu$) such that if $C=C(\pi)$ and}
$$
\alpha\le \frac{b_1}{\log C},
$$
\textit{then}
$$
\Big| \int_{-\infty}^{+\infty} \frac{L(1/2+it,\pi)}{1/2+it} g_\alpha(t)dt\Big| \ge c.
$$

\subsection{} We can now apply the previous proof by using a theorem of Molteni \cite{Mo}.

\vspace{2mm}

{\montheo (Molteni). Assume $\pi$ is a cuspidal, unitary representation of $GL(m,\A)$. Then}
$$
\sum_{n\le x} |a_n| = O_\varepsilon(C^\varepsilon x^{1+\varepsilon}) \qquad (x\ge 1).
$$

\noindent (The implicit constant depends only on $m$ and $\varepsilon$.) In fact Molteni proves a stronger result:
$$
\sum_{n\le x} \frac{|a_n|}{n} = O(C^\varepsilon x^{\varepsilon}).
$$
Thus, in the previous proof, we may take $\nu=1$, $\xi=\varepsilon$. (This has been improved by Xiannan Li \cite{Li}; the improved estimate seems irrelevant to us, but the theorem extends to automorphic representations satisfying the conditions in \S 1.1) The assumption $\sigma_c(\pi)\le \nu$ is of course satisfied.

We have given an exposition of the proof valid for other exponents, because of the following fact. Consider only representations $\pi$ that are tempered, i.e. satisfy  the Ramanujan Conjecture. Assume moreover $\pi_\infty$ self--dual. Then by a result of Friedlander and Iwanie\v c, (2.1) is satisfied with $\xi=\frac{1}{m+1}$, $\nu=\frac{m-1}{m+1}$. (See \cite{FI}; the self--duality condition is implicit there, and made explicit in \cite{Cl}.) However Friedlander and Iwanie\v c consider only the arithmetic conductor, so (2.1) is replaced by
$$
A_0(x) \ll D^{\frac{1}{m+1}} x^{\frac{m-1}{m+1}+\varepsilon}.
$$
They must also assume that $\pi_{\infty}$ is bounded, i.e. $C_{\infty}(\pi)=\prod (2+ \vert c_i \vert )$ bounded. Their estimate with respect to $x$ is better, but this does not seem to play any role in the present proof.

The fact that these  representations $\pi$ verify $\sigma_c(\pi) \le \frac{m-1}{m+1}$ is proved in \cite{Cl}.

We will unfortunately have to assume the Ramanujan conjecture in the next paragraph.

\subsection{} In order to exploit Proposition 2.1 to obtain a lower bound on a ``short'' interval, we must now estimate the tail
$$
I_5(T) = \int_T^\infty \frac{L(1/2+it)}{1/2+it} g_\alpha(t)dt
$$
(and the opposite one). In order to have uniform estimates, we must now assume that $\pi$ verifies the Ramanujan Conjecture. We now have
$$
L(\frac{1}{2}+it) \ll C(t)^{\frac{1}{4}+\varepsilon}.
$$
Fix $g_\alpha(t)= e^{-\pi\alpha t^2}$ with $\alpha=\frac{b_1}{\log C}$. For $T\ge 1$,
\begin{equation}
I_5(T) \ll C^{\frac{1}{4}+\varepsilon} \int_T^\infty e^{-\pi\alpha t^2} t^{\frac{m}{4}-1+\varepsilon}dt.
\end{equation}

We need an estimate for
$$
\int_T^\infty e^{-\pi\alpha t^2}t^a dt;
$$
as recalled in Chapter 1,
$$
\begin{array}{c}
\displaystyle \int_T^\infty e^{-\pi\alpha t^2} t^adt = \frac{1}{2} \pi^{-1} \alpha^{-1} T^{a-1} e^{-\pi\alpha T^2}\cdot \\
\cdot \displaystyle \Big( 1+O\Big(\frac{1}{\pi\alpha} T^{-2}\Big)\Big).
\end{array}
$$

We want to ensure that $|I_5(T)| \le \frac{1}{4}c$; i.e.
$$
C^{\frac{1}{4}+\varepsilon} \alpha^{-1} T^{a-1} e^{-\pi\alpha T^2} \Big(1+O\Big(\frac{1}{\alpha} T^{-2}\Big)\Big) \ll 1.
$$
We assume provisionally that $T$ can be chosen  such that $\frac{1}{\alpha}T^{-2}$ is small. Then this is implied~by
\begin{equation}
\frac{1}{2}\log C-\log \alpha +(a-1)\log T -\pi\alpha T^2 \le -X
\end{equation}
where $X$ is a positive constant; here $a=\frac{m}{4}-1+\varepsilon$; we replace $a$ by a value such that $a-1>0$. With $\alpha=\frac{b_1}{\log C}$, we want to ensure
$$
\pi\alpha T^2 -a \log T \ge X-\log\alpha +\frac{1}{2}\log C,
$$
which yelds
$$
\pi b_1T^2 -a\log C \log T \ge \log C (X'+\log\log C+\frac{1}{2}\log C).
$$
Assume $T^2=A_1\log^2 C$ for a large positive constant. The left--hand side is then
$$
\pi b_1 A_1 \log^2 C - a\log C(\log\log C+\log \sqrt{A_1}) \ge A_2 \log^2 C
$$
for $A_1$ sufficiently large and $C\ge 2^m$, and this dominates the right--hand side. Moreover, $\frac{1}{\alpha}T^{-2} \ll (\log C)^{-1} A_1^{-2} \le \frac{1}{2}$ for $A_1$ sufficiently large (independently of $C$), as assumed previously. Therefore:

\vspace{2mm}

{\monlem For $T=A\log C$, where $A$ is a sufficiently large positive constant, and $\alpha=\frac{b_1}{\log C}$,}
$$
\Big|\int_{|t|\ge T} \frac{L(1/2+it)}{1/2+it} g_\alpha(t)dt\Big| \le \frac{c}{2}.
$$

Now the absolute value of the integral for $|t|\le T$ is $\ge \frac{c}{2}$. Thus
\begin{equation}
\overset{}{\underset{t\in[-T,T]}{\Sup}} \Big|\frac{L(1/2+it)}{1/2+it}\Big| \ge \frac{c}{2}\ \Big|\int_{t\le T} g_\alpha(t)dt\Big|^{-1}
\end{equation}
\begin{equation}
\Big\|\frac{L(1/2+it)}{1/2+it}\Big\|_2  \ \Big\| g_\alpha(t)\Big\|_2 \ge \frac{c}{2}
\end{equation}
where the $L^2$--norms are computed in $[-T,T]$.

We now have to estimate the integral, and quadratic integral, of $g_\alpha(t)$ on $[-T,T]$. We have
$$
\begin{array}{rl}
\displaystyle \int_0^T e^{-\pi\alpha t^2} dt &=\displaystyle\int_0^{\sqrt{\pi\alpha}T} e^{-s^2}(\sqrt{\pi\alpha})^{-1} ds\\
&\displaystyle =a_6(\log C)^{1/2} \int_0^{a_7(\log C)^{1/2}} e^{-s^2}ds\\
&\asymp (\log C)^{1/2}
\end{array}
$$
and the same estimate is true for the quadratic integral. Therefore we obtain

\vspace{2mm}

{\montheo There exist absolute, positive constants $A$, $c_1$, $c_2$ such that for any cuspidal $\pi$ verifying the Ramanujan Conjecture

\vspace{2mm}

\noindent (i) $\overset{}{\underset{t\in[-T,T]}{\Sup}} \Big|\dfrac{L(1/2+it),\pi}{1/2+it}\Big| \ge c_1(\log C)^{-1/2}$
\vspace{2mm}

\noindent(ii) $\displaystyle \int_{-T}^T \Big|\dfrac{L(1/2+it),\pi}{1/2+it}\Big|^2 dt \ge c_2(\log C)^{-1/2}$

\vspace{2mm}

 if $T\ge A\log C$, $C=C(\pi)$.}

\vspace{2mm}

\vspace{2mm}

It may be noticed that in $(i)$ the Lindel\"of Conjecture yields an estimate $\ll (\log C)^\varepsilon $ (for $T$ of order $A \log C$).

We have used the cuspidality of $\pi$ only in  order to rely on Molteni's result; Xianan Li's result makes this  unnecessary.

We can now use the remarks in \S 1.3 on the translation of integrals to obtain lower bounds on $L(\frac{1}{2}+it,\pi)$ on short intervals, not necessarily centered in $0$. We have
$$
\begin{array}{c}
L(s,\pi[X] = L(s+i X,\pi)\\
\noalign{\vskip1mm}
C(\pi[X]) \le C(\pi) (1+|X|)^m\\
\noalign{\vskip1mm}
\log C(\pi[X]) \le \log C(\pi)+m \log (1+|X|).
\end{array}
$$
Let
$$
S= \overset{}{\underset{t\in[-T,T]}{\Sup}} \Big|\frac{L(1/2+it,\pi[X])}{1/2+it}\Big|.
$$

By Theorem 2.2, for $T\ge A \log C(\pi[X])$:
$$
S\ge c_1(\log C+ m\log (1+|X|))^{-1/2}.
$$
However
$$
S=\overset{}{\underset{t\in[-T,T]}{\Sup}} \Big|\frac{L(1/2+i(t+X),\pi}{1/2+it}\Big| \le 2 \overset{}{\underset{t\in[-T,T]}{\Sup}} |L(1/2+i(t+X,\pi)|
$$
Thus we have the following result on ``short intervals''.

{\montheo There exist absolute, positive constants $A$, $c_3, c_4 >0$ such that for any cuspidal $\pi$ verifying the Ramanujan Conjecture

\vspace{2mm}

\noindent(i)  $ {\underset{t\in[X-T,X+T]}{\Sup}} |L(1/2+it,\pi)|\ge c_3\cdot (\log C(\pi) +m\log (1+|X|))^{-1/2}.$

\vspace{2mm}
\noindent(ii)  $\displaystyle \int_{X-T}^{X+T} \Big|\dfrac{L(1/2+it,\pi}{1/2+i(t-X)}\Big|^2 dt \ge c_4 (\log C(\pi) +m\log (1+|X|))^{-1/2}$

\vspace{2mm}

if $T\ge A(\log C+ \log (1+|X|)$, $C=C(\pi)$.}

\vspace{2mm}

The absolute constants depend only on $m$. 

Sarnak suggested that the estimate (ii) should be compared with the lower bound on short intervals obtained by Ramachandra in \cite{Ra}. (Ramachandra notes that his proof, written for the Riemann zeta function, will extend to $L$-functions with Euler products. See \cite[Remark 3]{Ra}.) In fact, because of the large denominator in the integral, our lower bound seems better than the one that follows obviously from Ramachandra's. We assume $X$ positive and large and $T$ small with respect to $X$. (In the next lines the constant $A$ depends on the formula.) Then Ramachandra's estimate is

$$
 \displaystyle \int_{X}^{X+T} |L(1/2+it)|^2 dt \ge A T \log T, \\\\\  T \geq A \log \log X
 $$
 
 \vspace{2mm}
 
 \noindent while here we obtain
 
 $$
  \displaystyle \int_{X-T}^{X+T} \Big|\dfrac{L(1/2+it,\pi}{1/2+i(t-X)(}\Big|^2 dt \ge A (\log C + \log X)^{-1/2}
$$

 \noindent for $T \geq A(\log C+ \log X).$ 

For fixed $C$, our condition on $T$ implies Ramachandra's.  We consider $T \approx A\log X$. The denominator in our formula  is smaller than $T^2$. Thus Ramachandra's bound implies for the  integral

 $$
  \displaystyle \int_{X-T}^{X+T} \Big|\dfrac{L(1/2+it,\pi}{1/2+i(t-X)}\Big|^2 dt 
$$

 \noindent a lower bound in $A T^{-1} \log T$. But for $T \approx A\log X$, $(\log X )^{-1} \log \log X$ is dominated by our bound $(\log X)^{-1/2}$. Thus, for fixed $C$, the results do not seem commensurable.  (It would be interesting to know how the conductor enters in Ramachandra's formula.)

\section{The first non--trivial Frobenius element in a Galois representation}

\subsection{} 
In this chapter we consider a non--trivial, irreducible representation $\rho:\Gal(E/F) \rg GL(m,\C)$, where $E/F$ is a Galois extension of number fields. If $\rho$ is unramified at a prime $\pp$ of $F$, we can consider the image $\rho(\Frob_\pp)$, a conjugacy class in $GL(m,\C)$. We want to obtain a lower bound $\beta$ on the first value of $N\pp$ such that $\rho$ is unramified and $\rho(\Frob_\pp)\not=1$. Of course, this is complicated by the presence of ramification, i.e., there may be a ramified prime $\fq$ such that $N\fq\le\beta$. If $m=1$, this applies to Artin characters of~$\A_F^\times$.

\subsection{}

We first consider the Abelian case, $m=1$. If $F=\Q$, we are looking at primitive Dirichlet characters $\mod\, q=D$. If $L(s,\chi)\ll (q|s|)^{\mu+\varepsilon}$ for $Re(s)=\frac{1}{2}$, it was shown in \cite[\S~3.3]{Cl} that $\beta\ll q^{2\mu+\varepsilon}$. Since Petrow  and Young \cite{PY} have proved the Weyl bound $\mu\le\frac{1}{6}$, we obtain by this simple application of the Mellin transform
\begin{equation}
\mathrm{For}\ \chi\ \mathrm{a\ primitive\ Dirichlet\ character} \mod q, \ \beta(x) \overset{}{\underset{\varepsilon}{\ll}} q^{1/3+\varepsilon}.
\end{equation}
Of course this is much weaker than Burgess's bound (for $q$ prime)
$$
\beta(\chi)\overset{}{\underset{\varepsilon}{\ll}} q^{\frac{1}{4\sqrt{e}}+\varepsilon}.
$$
However, we can now extend the method to Artin characters of $\A_F^\times$ for an arbitrary number field~$F$.

So let $\chi$ be a character of finite order, and $\rho$ the $1$--dimensional representation of $\Gal(F_{\chi}/F)$ associated to $\chi$ by class field theory. We consider the $L$--function $L(s,\chi)$ as a $L$--function over $\Q$. All the conditions in the paper of Friedlander--Iwaniec are met --- this will also be the case for a Galois representation of degree $>1$. We therefore have  the properties of $L(s,\chi)$ recalled in \S~2.3.

Rather than the formulas in \cite{Cl}, we use the Mellin transform as in \S~2.1. Thus let
$$
\begin{array}{rl}
L(s,\chi) &=\displaystyle\sum_n a_nn^{-s},\\
\noalign{\vskip2mm}
A_0(x) &=\displaystyle \sum_{n\le x}a_n
\end{array}
$$
and
\begin{equation}
\widetilde{\Mm} \, A_0(x)=\int_1^\infty A_0(x) x^{-s}\frac{dx}{x}.
\end{equation}

Let $\nu_0$ be the exponent of Friedlander--Iwaniec:
$$
\nu_0=\frac{d-1}{d+1}, \ d=[F:\Q].
$$
so
\begin{equation}
A_0(x) {\underset{d, \varepsilon}{\ll}} D^{\frac{1}{d+1}} x^{\nu_0+\varepsilon}\quad (x\ge 1).
\end{equation}

We note that the Archimedean data $(c_j)$ associated to $\chi$ are uniformly bounded, so the constant depends only on $(d,\varepsilon)$. Moreover (3.2) is obtained by Friedlander--Iwaniec for $x\ge D^{1/2}$ (see a correction in \cite[\S~2.5]{Cl}) but it is trivial for $x\le D^{1/2}$ since $|a_n|\ll n^{\varepsilon'}$ for any $\varepsilon'>o$.

The integral in (3.2) is absolutely convergent for $Re(s)>\nu_0$, say $Re(s)=\nu>\nu_0$. The function $A_0(x) x^{-\nu}$ then belongs to $L^2(\R_+^\times,\frac{dx}{x})$. It follows (cf. Titchmarsh \cite{Ti}) that $\frac{L(s)}{s}$ is $L^2$ on the line $Re(s)=\nu$, and that
\begin{equation}
\int_1^\infty |x^{-\nu} A_0(x)|^2 \frac{dx}{x} = \frac{1}{2\pi} \int_\nu \Big|\frac{L(s)}{s}\Big|^2\ |ds|.
\end{equation}
However, if $\nu$ is close to $\nu_0$, the convergence of the right--hand side does not follow from convexity: see \cite[\S~3.2]{Cl}; we will review this in Chapter~4. As a consequence we cannot use Proposition~1.1 to estimate the right--hand side.

By Proposition 1.1, we have $\mu(\nu)\le \frac{d}{2}(1-\nu)$, so $\mu(\nu)<\frac{1}{2}$ if $\nu> 1-\frac{1}{d}$. The convergence of the right--hand side then follows from convexity and we have again by Proposition~1.1:

\vspace{2mm}

{\monlem For $\nu>1-\frac{1}{d}$,}
$$
\int_\nu \Big|\frac{L(s)}{s}\Big|^2 \, |ds| \ll D^{1-\nu+\varepsilon}.
$$

\subsection{}
 
 We now want to obtain a lower bound for the left--hand side of (3.4). We consider the integral on $[1,\beta]$ where $\beta$ is the norm of the first unramified prime such that $\chi(\pp)\not=1$.
 
 Before doing so we recall the expression of the ``absolute" conductor $D$. Let $D_F$ be the absolute value of the discriminant of $F$; let $\ff(\chi)$ be the conductor of $\chi$ (seen as a character of $\A_F^\times$), an ideal of $F$. Then
 \begin{equation}
 D=D_F\ N_{F/\Q} \ff(\chi).
 \end{equation}
 
 We need an expression for $A_0(x)$, for $x<\beta$. It is then given by
 $$
 a_n =\sum_{N\fa=n} a(\fa)
 $$
 where $\sum a(\fa)N\fa^{-s}$ is the expression of $L(s,\chi)$ as an $L$--function over $F$, and $N\fa < \beta$, which implies that the factorization of $\fa$ involves only primes $\pp$ with $N\pp<\beta$. We have
 $$
 L_\pp (s,\chi)=(1-N\pp^{-s})^{-1} 
$$
if $\chi$ is unramified at $\pp$ and $N\pp<\beta$,
$$
L_\fq(s,\chi)=1
$$
is $\chi$ is ramified at $\fq$. If $\fa=\prod \fq_i^{\alpha_i} \prod \pp_j^{\beta_j}
$
(where we use $\{\fq_{i}\}$ to denote all the ramified primes),  this implies that
$$
a(\fa)=0
$$
if $\alpha_i>0$ for some $i$, and $a(\fa)=1$ otherwise. Thus, up to $x=\beta$, $A_0(x)$ is a variant of the summation function for $\zeta_F$:
\begin{equation}
A_0(x) = \sum_{n\le x} \sum_{N\fa=n\atop (\fq_i,\fa)=1} 1
\end{equation}
where the condition $(\fq_i,\fa)=1$ is imposed for all $i$. Let $B_0(x)$ be the summation function for $\zeta_F$. The difference $B_0-A_0$ is then
\begin{equation}
 \sum_{n\le x} \sum_{N\fa=n}{}' 1
\end{equation}
where $\sum'$ is restricted by the condition
\begin{equation}
\exists ~~ i : \fq_i|\fa.
\end{equation}

At this point it is necessary to extend the main result of \cite{FI}. Very generally, assume $\pi$ is a tempered representation of $GL(m,\A)$; assume $\pi$ is a product $\pi_1\times\cdots\pi_r$ of cuspidal representations (with at most a factor $m_i=1$, $\pi_i=1$); assume $\pi_\infty$ self-dual. Let $A_0(x)$ be the associated summation function. Let $\kappa$ be the residue of $L(s,\pi)$ at~1. We write $\gD$ for the conductor of $\pi$.

\vspace{2mm}

{\maprop For $x\ge1$,
$$
A_0(x) = \kappa x+O(\gD^{\frac{1}{m+1}+\varepsilon} x^{\frac{m-1}{m+1}+\varepsilon})
$$
where the implicit constant depends only on $m$, $\varepsilon$ and the Archimedean parameters of~$\pi$.}

\vspace{2mm}

For the function $\zeta_F$, the implicit constant depends only on $d$ and~$\varepsilon$.

For $x\ge \gD^{1/2}$ this is the result of \cite{FI}.  For $x\le \gD^{1/2}$, it was already noticed in \S~3.2 that $A_0(x) \ll \gD^{\frac{1}{m+1}} x^{\frac{m-1}{m+1}+\varepsilon}$. (This does not depend on the fact that $L(s)$ is holomorphic.) If $x\le \gD^{1/2}$, $\gD^{\frac{1}{m+1}}x^{-\frac{2}{m+1}}\ge 1$ and it suffices to check: 
$$
\kappa\ll (\gD x)^\varepsilon
$$
which is true since $\kappa\ll \gD^\varepsilon$, cf. \cite[p. 100]{IK}. (For a z\^eta function, this is the easy part of the Brauer--Siegel theorem.)

\vspace{2mm} 

We can now apply this to $B_0$:
$$
B_0(x) = \kappa x+O(D_F^{\frac{1}{d+1}+\varepsilon} x^{\frac{d-1}{d+1}+\varepsilon}),
$$
$\kappa=Res_{s=1}\zeta_F(s)$. Now fix $i$ and consider the sum (3.7), the condition being $\fq_i|\fa$. By a change of variables,
$$
\sum_{N\fa\le x\atop \fq_i|\fa} 1 = \sum_{N\fb \le \frac{x}{q_i}}1
$$
where $q_i=N\fq_i$ and $\fb$ ranges over integral ideals. This can be expressed as
$$
\kappa \frac{x}{q_i} +O(D_F^{\frac{1}{d+1}+\varepsilon}(x/q_i)^{\frac{d-1}{d+1}+\varepsilon}).
$$
However, care has to be exercised. This is certainly true for $x\ge q_i$, however for $x<q_i$ the sum is empty and we have to check that
$$
\kappa\frac{x}{q_i} \ll D_F^{\frac{1}{d+1}+\varepsilon}(x/q_i)^{1-\frac{2}{d+1}+\varepsilon}
$$
(uniformly when $d$ is fixed.) This reduces to
$$
\kappa \ll D_F^{\frac{1}{d+1}+\varepsilon}(x/q_i)^{-\frac{2}{d+1}+\varepsilon};
$$
but $x/q_i\le 1$ and $\kappa\ll D_F^\varepsilon$.

Now return to the sum (3.7). By the  inclusion--exclusion principle, it is equal to
$$
\sum S_i - \sum S_{ij}+ \sum S_{ijk} - \cdots
$$
where
$$
S_{ij\cdots \ell} = \sum_{N\fa\le x\atop \fq_i\cdots \fq_\ell|\fa} 1.
$$
The same argument yields, with $\fq=\fq_i\cdots \fq_\ell$, $q=q_i\cdots q_\ell$:
$$
S_{ij\cdots \ell} = \kappa \frac{x}{q_i} +O\Big(D_F^{\frac{1}{d+1}+\varepsilon}(x/q)^{\frac{d-1}{d+1}+\varepsilon}\Big).
$$
The previous argument remains correct for $x\le q$, with the same uniformity.

Let us write $\xi=\frac{d-1}{d+1}+\varepsilon$. We now have
\begin{equation}
A_0(x) = \prod_i\Big(1-\frac{1}{q_i}\Big) \kappa x+ O(\prod (1+q_i^{-\xi})D_F^{\frac{1}{d+1}+\varepsilon}x^\xi).
\end{equation}
In (3.9), we can (brutally) replace $D_F$ by $D\geq D_F$. 
We have to evaluate $Q=\prod(1-\frac{1}{q_i})$ and $R=\prod(1+q_i^{-\xi})$.

Let $p_i$ be the prime divisor of $q_i$. Then
$$
Q\ge \prod_i \Big(1-\frac{1}{p_i}\Big).
$$
There may be several primes associated to $p_i$, but fewer than $d=[F:\Q]$. Thus
$$
Q\ge \prod_j \Big(1-\frac{1}{p_j}\Big)^d
$$
where the $p_j$ are now distinct. But $p_j$ divides $N \ff_\chi$, so $p_j\le N=N \ff_\chi$ and
$$
Q\ge \prod_{p\le N} \Big(1-\frac{1}{p}\Big)^d.
$$
Since $\prod\limits_{p\le N}(1-\frac{1}{p}) = \frac{e^{-\gamma}}{\log N}(1+O(\frac{1}{\log N}))$ we deduce that $Q\ge A\, \log(N \ff_\chi)^{-d}$. (If $\chi$ is everywhere unramified, $Q=1$.).

Consider now  $R=\prod (1+q_i^{-\xi})$. Again, $R\le \prod\limits_i (1+p_i^{-\xi}) \le (\prod\limits_j (1+p_j^{-\xi}))^d$, with the same notation. But $\prod(1+p_j^{-\xi})=\sigma_{-\xi}(p_1\cdots p_j) < \sigma_0(p_{1}\cdots p_j) \ll (p_1\cdots p_j)^\varepsilon \ll (N \ff_\chi)^\varepsilon$ and the same estimate obtains for~$R$.

Now fix $\nu>1 -\frac{1}{d}$ and consider the integral
$$
\int_1^\beta |x^{-\nu}A_0(x)|^2\frac{dx}{x}.
$$
By (3.9) this is a sum of three terms, one of them exact:
$$
I_1= \int_1^\beta(Q \kappa x^{1-\nu})^2 \frac{dx}{x},
$$
so
\begin{equation}
I_1\asymp Q^2 \kappa^2\beta^{2-2\nu} \qquad (\beta\ge 2).
\end{equation}

The second term is dominated by
$$
I_2 = D^{\frac{1}{d+1}+\varepsilon} Q\kappa \int_1^\beta x^{2-\frac{2}{d+1}-2\nu} \frac{dx}{x}
$$
using that $R\ll D^\varepsilon$, so
\begin{equation}
I_2 \asymp D^{\frac{1}{d+1}+\varepsilon} Q \kappa \beta^{2-\frac{2}{d+1}-2\nu}.
\end{equation}

The third term is dominated by
$$
I_3 = D^{\frac{2}{d+1}+\varepsilon} \int_1^\beta x^{-2\nu+1-\frac{4}{d+1}+\varepsilon} dx
$$
(rescaling $\varepsilon\cdots$).

For $\nu>1-\frac{1}{d}$, the exponent in the integral is always $<-1$ (Assume $d>1$.)  Thus
\begin{equation}
I_3 \asymp D^{\frac{2}{d+1}+\varepsilon}.
\end{equation}
We now compare $I_1$ and $I_2$; we want to assume that $I_1$ is dominant, i.e.
$$
Q\kappa \succ D^{\frac{1}{d+1}+\varepsilon} \beta^{-\frac{2}{d+1}}.
$$

\textit{We now assume that the extension $F/\Q$ is Galois.} By the Brauer--Siegel theorem, $\kappa\gg D_F^{-\varepsilon}$, the implicit constant depending only on $d$. By our estimate on $Q$, we must have
$$
D_F^{-\varepsilon}(\log N \ff_\chi)^{-d} \succ D^{\frac{1}{d+1}+\varepsilon} \beta^{-\frac{2}{d+1}}
$$
which is true, upon changing $\varepsilon$, if
$$
\beta^{\frac{2}{d+1}} \succ D^{\frac{1}{d+1}+\varepsilon},
$$
i.e.
\begin{equation}
\beta
 \gg D^{\frac{1}{2}+\varepsilon}.
 \end{equation}

Now consider $I_1$ and $I_3$. In this case, using again the estimate on $Q\kappa$, we must have
$$
\beta^{2-2\nu} \succ  D^{\frac{2}{d+1}+\varepsilon};
$$
for $\nu$ close to $1-\frac{1}{d}$, we see that this yields
\begin{equation}
\beta
\gg D^{\frac{d}{d+1}+\varepsilon}.
\end{equation}

We now compare $I_1$, assumed to be dominant, with the estimate given by Lemma~3.1. This yields
$$
D^{1-\nu+\varepsilon} \gg Q^2 \kappa^2 \beta^{2-2\nu};
$$
since we can neglict  $Q$ $\kappa$, we obtain
$$
\beta \ll D^{1/2+\varepsilon}.
$$
The conclusion is:
$$
\beta \gg D^{\frac{d}{d+1}+\varepsilon} \Rightarrow \beta \ll D^{1/2+\varepsilon}.
$$
We conclude that $\beta\ll D^{\frac{d}{d+1}+\varepsilon}$ for large $D$. There is a finite number of pairs $(F, \chi)$ such that $D_F N(\ff_\chi)\le A$. Indeed there is then a finite number of possibilities for $F$; for $F$ fixed the primes and the ramification degrees of $\chi$ are bounded. Therefore:

\vspace{2mm}

{\montheo Consider the Galois extensions $F$ of degree $d\ge 2$ of $\Q$, and the non--trivial Artin characters $\chi$ of $F$. Let $\beta=N\pp$ where $\pp$ is an unramified prime of smallest norm such that $\chi(\pp)\not=1$. Then
$$
\beta \ll D^{\frac{d}{d+1}+\varepsilon}
$$
where $D=D_F N(\ff_\chi)$ and the implicit constant depends only on $d$ and $\varepsilon$.}

\vskip2mm

\noindent \textbf{Remark. } It may be possible to obtain a slightly better estimate, depending of $N \ff_\chi$, by avoiding the change of $D_F$ to $D$ after (3.9).

\subsection{} Now let $F$ be a number field and
$$
\rho : \Gal(E/F) \lgr GL(M,\C)
$$
an irreducible representation of degree $M$. We consider its Artin $L$--function $L(s,\rho)$, given by an Euler product of degree $M$ over $F$. We can view it as an $L$--function of degree $m=M d$ over $\Q$. We assume that $\rho$ is non--trivial and that the Artin Conjecture is true for $\rho:L(s,\rho)$ is holomorphic. (For a review of recent results on the Artin Conjecture see Calegari [Ca ] and the references therein.). We recall that the results of Friedlander--Iwaniec apply to $L(s,\rho)$ (cf. \cite[\S~2.1]{Cl}.)

The (absolute) conductor $D$ of $L(s,\rho)$ is equal to $D_F^M\, N(\ff_\rho)$ where $\ff_\rho$, an ideal of $F$, is the Artin conductor of $\rho$, cf. Neukirch~\cite[Prop.11/7]{Neu}.

Let $L(s,\rho)=\sum\limits_1^\infty a_n n^{-s}$. Its Euler product (over $F$) is 
$$
L(s,\rho) = \prod_\pp \det((1-\Frob_\pp N\pp^{-s})\mid V^{I_\pp})^{-1}
$$
where $V$ is the space of $\rho$ and $I_\pp$ is the inertia. In particular, the ramified primes introduce factors $\prod\limits_i (1-\alpha_{i,\pp}N \pp^{-s})^{-1}$, where $\alpha_{i,\pp}$ is a root of unity different from 1, that are not positive. In order to apply the method of \S~3.2, we introduce the \textit{unramified $L$--function}
$$
\Ll(s,\rho) = \prod_\pp L_\pp(s,\rho)
$$
where the product ranges only over the unramified  primes. Thus, with $\fq$ ranging over the ramified primes,
$$
\Ll(s,\rho) = L(s,\rho) \prod_\fq D_\fq (s,\rho)
$$

where $D_\fq(s,\rho)=\prod\limits_i (1-\alpha_{i,\fq} N\fq^{-s})$. Let $D(s,\rho) =\prod\limits_\fq D_\fq(s,\rho)$.
We need to control the growth of $\Ll(s,\rho)$ in the critical strip. For\break $Re(s)\ge 0$, $D_\fq(s,\rho) \le 2^{M-1}$. Thus
$$
|D(s,\rho)| \le 2^{(M-1)r} \le 2^{d(M-1)t}
$$
where $r$ is the number of ramified primes $\fq$ and $\{p_1,\dots p_t\}$ are the distinct rational primes dividing one of the $\fq_j$.

Fix $d$, and fix $\varepsilon>0$. We first show that
\begin{equation}
|D(s,\rho)| \le (D_F^M \, N \ff_\rho)^\varepsilon
\end{equation}
except for a finite number of pairs $(F,\ff)$. Assume first $r\rg \infty$. We have $r\le dt$, so $t\ge [\frac{r}{d}]$. Moreover $p_i\mid N\ff_\rho$, so $N\ff_\rho \ge \prod p_i\ge t$! Thus (3.14) is verified if $2^{(M-1)r} \le \Gamma([\frac{r}{d}]+1)^\varepsilon$, thus if
$$
2^{M-1)r} \le \Big(\Gamma\Big(\frac{r}{d}\Big)\Big)^\varepsilon.
$$
Since $\varepsilon$ is fixed, this is true for large $r$ by Stirling's formula.

So assume now $r\le r_0$. Then  $|D(s,\rho)|\le 2^{(M-1)r_0}$ and (3.15) is true if $D_F$ is large, i.e., excluding a finite number~of~$F$.

Again, for $F$ fixed, (3.15) is true if $N\ff_\rho$ is sufficiently large. If $N\ff_\rho\le A$, this leaves a finite number of possibilities for the places and degrees of~$\ff_\rho$.

Finally (3.15) is violated only for a finite number of values of $(F,\ff_\rho)$. Therefore we have
$$
|D(s,\rho)| \le A(D_F^M N\ff_\rho)^\varepsilon
$$
for a sufficiently large constant $A$. By  Proposition~1.1, then:

\vspace{2mm}

{\monlem For $Re(s)\ge 0$,
$$
\Ll(s,\rho) \ll C(s)^{\frac{1-\sigma}{2}+\varepsilon}
$$
where $C(s)$ is the analytic conductor of $L(s,\rho)$.}

\vspace{2mm}

We can now imitate the arguments of \S 2. Let $\beta$ be the smallest norm $N\pp$ of an unramified prime such that $\rho(\Frob_\pp)\not=1$. If $\Ll(s,\rho) =\sum\limits_\fa b_1 (\fa) N\fa^{-s}$ is the Dirichlet series of $\Ll(s,\rho)$ over $F$, the coefficients $b_1(\fa)$ coincide, for $N\fa<\beta$, with  those of the Euler product $(S=\{\fq_1,\ldots, \fq_s\})$
$$
\begin{array}{rl}
L^S(s,\rho)& = \prod\limits_{\pp\notin S} (1-N\pp^{-s})^{-M}\\
&=(\zeta_F^S(s,\rho))^M.
\end{array}
$$
These coefficients are equal to $0$ if $\fa$ is divisible by one of the $\fq_i$, and $\ge1$ otherwise. (Using the explicit series for $(1-X)^{-M}$ does not lead to better estimates.). Thus, with $\fq=\prod \fq_j$:
\begin{equation}
\Ll(s,\rho) = \sum_{n=1}^\infty b_n \, n^{-s}
\end{equation}
where
\begin{equation}
b_n \ge \sum_{N\fa=n\atop (\fa,\fq)=1} 1\quad (n<\beta).
\end{equation}

Now let $B_0$ be the summation function of $\Ll(s,\rho)$:
$$
B_0(x) = \sum_{n\le x}b_n.
$$

For $x<\beta$, this is given by (3.16), and this has been analysed in \S~3.3. (See formulas (3.6) to (3.9).) In fact, with $A_0(x)$ defined by (3.6), we now have

\vspace{2mm}

{\monlem For $x<\beta$,
$$
B_0(x) \ge A_0(x)
$$
where (for $x<\beta$)}
$$
A_0(x) =Q \kappa x +O(R D^{\frac{1}{d+1}+\varepsilon}x^{\frac{d-1}{d+1}+\varepsilon}).
$$

\vspace{2mm}

Again, we have replaced $D_F$ by $D$; $\kappa$ is the residue of $\zeta_F$.

However, we still have to check that the Mellin transform can be applied  to yield the summation function $B_0(x)$, and for which values of~$s$.

\subsection{}

 We write simply $D(s,\rho)=\prod\limits_j (1-\alpha_j N\fq_j^{-s})$ where $\fq_j$ ranges over ramified primes and $\alpha_j$ over the associated roots. Let $q_j=N\fq_j$. With $L(s,\rho)=\sum\limits_n a(n)n^{-s}$,
 $$
 L(s,\rho) (1-\alpha q^{-s}) = \sum_{n\ge 1} (a(n)-\alpha a(n/q))n^{-s}
 $$
where $\alpha_{n/q}=0$ if $q \nmid n$. We deduce that $\Ll(s,\rho)= \sum\limits_n b_n n^{-s}$, where
$$
\begin{array}{c}
b_n =a(n) -\sum\limits_j \alpha_j a(n/q_j)+\sum\limits_{j_1,j_2} \alpha_{j_1}\alpha_{j_2} a(n/q_{j_1}q_{j^2}) -\\
\cdots +(-1)^N \alpha_1 \cdots \alpha_N a(n/q_1\cdots q_N),
\end{array}
$$
$N$ being the degree of $D(s,\rho)$; the same condition on  $a(n/q_1\cdots q_n)$ applies.

Recall that $m=dM$. By the theorem of Friedlander--Iwaniec, we have 
$$
\ga_0(x) \ll D^{\frac{1}{m+1}+\varepsilon} x^{\frac{m-1}{m+1}+\varepsilon} \qquad (x\ge 1)
$$
where $\ga_0$ is the summation function for $L(s,\rho)$,
and we deduce that the same estimate is true for $B_0(x)$ if $x$ is sufficiently large. In particular (Lemma~2.1) $\widetilde{\Mm} B_0(s)$ is given by an absolutely convergent integral if $Re(s) >1-\frac{2}{m+1}$, and, for $\nu=Re(s)>1 - \frac{1}{m}$,
\begin{equation}
\int_\nu \Big|\frac{\Ll(s)}{s}\Big|^2 |ds| \ll D^{1-\nu+\varepsilon}
\end{equation}
(See Lemma 3.1.)

By Lemma 3.3, the integral
$$
\int_1^\beta |x^{-\nu}B_0(x)|^2 \frac{dx}{x}
$$
is larger than
$$
\int_1^\beta |x^{-\nu}A_0(x)|^2 \frac{dx}{x}
$$
where the expression of $A_0$ is recalled in Lemma~3.3. We now imitate the calculation in \S~3.3. We find first an explicit term, cf~(3.10)
$$
I_1 \asymp Q^2 \kappa^2\beta^{2-2\nu}.
$$

However, in the computation of $I_2$, the exponent $2-\frac{2}{d+1}-2\nu$ is negative for $\nu>1-\frac{1}{m}$, so
$$
I_2 \ll D^{\frac{1}{d+1}+\varepsilon} Q\kappa R.
$$
The same applies to the third term,
$$
I_3 \ll D^{\frac{2}{d+1}+\varepsilon} R^2 .
$$

Note that $Q$, $R$ and $\kappa$ depend on $F$, not on $\rho$. However $D_F\ll D$, so we can use the estimates in \S~3.3, and neglect these terms if $F$ is Galois over $\Q$. With $\nu \sim 1-\frac{1}{m}$, we see that $I_1$ dominates $I_2$ if 
\begin{equation}
\beta
 \gg D^{\frac{m}{2(d+1)}+\varepsilon},
 \end{equation}
and that it dominates $I_3$ if
\begin{equation}
\beta
\gg D^{\frac{m}{d+1}+\varepsilon}.
\end{equation}
Using (3.18) we see that
$$
\beta\gg D^{\frac{m}{d+1}+\varepsilon} \Rightarrow \beta\ll D^{1/2+\varepsilon}.
$$
The conclusion, as in \S 3.3, is that
$$
\beta \ll D^{\frac{m}{d+1}+\varepsilon}.
$$

\vspace{2mm}

{\montheo Let $\rho:\Gal(E/F) \rg GL(M,\C)$ be an irreducible, non trivial representation and $\beta=N\pp$ be the smallest norm of a prime $\pp$ of $F$ such that $\rho(\Frob \pp)\not=1$. Then, if $L(s,\rho)$ is holomorphic and $F/\Q$ Galois, 
$$
\beta\ll D^{\frac{m}{d+1}+\varepsilon} 
$$
where $d=[F:\Q]$, $m=Md$, and $D=D_F^M N(\ff_\rho)$.}

\vspace{2mm}

The implicit constant depends only on $M$ and $d$, but it is not effective since we have used the Brauer--Siegel theorem.

\subsection{}

It is of course of interest to compare this result with earlier estimates. Assume $E$ is minimal, i.e., $\rho:\Gal(E/F)\rg GL(M,\C)$ is injective; let $G=\Gal(E/F)$ and $g=|G|$. If $\Cc$ is a conjugacy class in $G$, let $\CPp(\Cc)$ be the set of primes $\pp$ of $F$ such that $\rho$ is unramified at $\pp$ and $\rho(\Frob_\fq)\in \Cc$; let $\CPp_1(\Cc)$ be the subset composed of primes such that $f(\pp/p)=1$; i.e. $N\pp=p$. Let $\beta(\Cc)$, $\beta_1(\Cc)$ denote the smallest norm of elements in $\CPp(\Cc)$, $\CPp_1(\Cc)$. Then the following results are known:
\begin{equation}
\beta
(\Cc) \ll (\log D_E)^2(\log\log D_E)^4
\end{equation}
(Lagarias and Odlyzko, 1975, under the generalized Riemann hypothesis)
\begin{equation}
  \beta_1(\Cc) \ll D_E^A  
\end{equation} 
(Lagarias, Montgomery and Odlyzko, 1979, unconditionaly) which has been improved by Zaman (2017) to
\begin{equation}
 \beta_1(\Cc) \ll D_E^{40}.
 \end{equation} 
 In the last two results, the implicit constant is effective. See \cite{LO,LMO,Z}. The last estimate has been improved to $\beta_1(\Cc) \leq D_E^{16}$ for large $D_E$, see \cite{Ka}.
 
 Clearly $\beta_1(\Cc) \ge \beta(\Cc)$ and $\beta\le \overset{}{\underset{\Cc}{\Inf}}~ \beta(\Cc)$, where $\Cc$ runs over non--trivial conjugacy classes. Since the estimates are uniform with respect to $\Cc$, (3.20)--(3.22) give an upper bound for $\beta$. However,
 $$
 D_E = D_F^g N_{F/\Q}(\fd_{E/D}),
 $$
 equal by the F\"uhrerdiskriminantenproduktformel to
 \begin{equation}
 D_F^g \prod_\rho (N_{F/\Q} \ff_\rho)^{\dim \rho},
 \end{equation}
 to be compared with
 \begin{equation}
 D^{\frac{m}{d+1}}= D_F^{M^2\frac{d}{d+1}} N(\ff_\rho)^{M\frac{d}{d+1}},
 \end{equation}
$M=\dim \rho$. We have $M^2<g$, so the first factor of (3.25) is negligible with respect to that of (3.24). Similarly, since $M=\dim \rho$, the second factor of (3.25) is smaller that that of (3.24) relative to our chosen $\rho$. In general the new bound is much better. However, our argument applies only to fields $F$, and representations $\rho$, of bounded degrees $d, M$.

A more pertinent comparison, however, is with Zaman's paper \cite{Z2}, which gives an estimate for $\beta$ (the smallest norm of a prime $\frak{p}$ of $F$ that does not split in $E$; Zaman also assumes that  $\frak{p}$ is of degree 1 over $F$. )  For simplicity we only consider the case where $F=\Q$.  In this case \cite[Cor.1.2]{Z2}, with $g=[E:\Q]$, Zaman's estimate is
$$
\beta \ll D_E^{\frac{1+\varepsilon}{4A(g-1)}}
$$
$$
\ll(\prod_{\rho} (N\ff_{\rho})^{\dim \rho})^{\frac{1+\varepsilon}{4A(g-1)}}
$$
with $A \approx 1$,  to be compared with (this paper)
$$
\beta \ll N(\ff_{\rho_0})^{M/2 + \varepsilon}, \\\ M=\dim\rho_0
$$
where $\rho_0$ is our chosen representation. The comparison therefore depends on the complexity of the Galois group $G$, in particular its number of large representations and their ramification. For larger fields $F$ \cite[Theorem1.1]{Z2}, this is multiplied by a term at least of order $\exp(Ag(\log D_F)^2)$ (here $A$ is again an implicit constant) and the product will likely be larger that the estimate of Theorem 3.2. Note also that Zaman does not assume the extensions to be Galois.

\vspace{2mm}

\textbf{Example:} Assume $\rho$ is associated to a normalised newform $\gf$ of weight 1 on $\Gamma_0(q)$. Then $\ff_\rho=(q) \subset \Z$, and Theorem~3.2 yields $\beta\ll q^{1+\varepsilon}$. 

\section{Complements to \cite{Cl} and remarks on a paper of Friedlander--Iwanie\v c}

\subsection{} The method of Friedlander--Iwanie\v c, and the proof in the Appendix of \cite{Cl}, use only the $L$--function of a (putative) representation $\pi$. In particular, they apply to the Rankin $L$--function $L(s,\pi_1\boxtimes \pi_2)$ of two cuspidal representations of $GL(m_1,\A)$ and $GL(m_2,\A)$. In particular, we have:
\begin{equation}
\begin{array}{c}
The\ lower\ bounds\ on\ the\ quadratic\ integrals\\
of\ L(s)\ in\ Theorems\ A-D \ of\  \cite{Cl} are\ true\\
if\ L(s)=L(s,\pi_1\boxtimes\pi_2),\ \pi_1,\pi_2\  cuspidal.
\end{array}
\end{equation}

Assume moreover $\pi_1$, $\pi_2$ tempered, and $\pi_{1,\infty},  \pi_{2,\infty}$ self--dual. Let
$$\begin{array}{c}
L(s,\pi_1 \boxtimes\pi_2) = \sum a_n n^{-s}\ \mathrm{and}\\
A(x) = \sum\limits_{n\le x} a_n n^{-s}.\ \mathrm{Then}
\end{array}
$$
\begin{equation}
A(x) = \kappa x+O(D^{\frac{1}{m+1}}x^{1-\frac{2}{m+1}+\varepsilon})
\end{equation}
\textit{where $m=m_1m_2$, and $\kappa$ is the usual residue, and $D$ is the conductor of $L(s,\pi_1\times \pi_2)$.}
 
The estimation  of $A_s(x)$ in \cite[Thm. 2.2]{Cl}, also applies; for $\pi_1 \ncong \tilde{\pi}_2  \vert ~\vert^{ia}$, this yields:
\begin{equation}
\begin{array}{c}
The\ abscissa\ of\ convergence\ of\ L(s,\pi_1\boxtimes\pi_2)\\
 satisfies\ \sigma_c \le 1-\frac{2}{m+1}.
\end{array}
\end{equation}

\subsection{} However, the results of Chapter 2 on ``short" integrals do not apply in general to Rankin $L$--functions. Indeed Molteni uses the properties of $L(s,\pi\boxtimes\tilde{\pi})$; for $\pi=\pi_1\boxtimes \pi_2$, this would assume the properties of the $L$--function of a quadruple tensor product.

\subsection{} In \cite[\S 3.4]{Cl} we discussed the relation between the exponent $\nu$ of an estimate of $A_0(x)$:
$$
A_0(x) =O(x^{\nu+\varepsilon}
)
$$
where $A_0(x)$ is associated to a cuspidal $\pi$, and the convexity estimates. For simplicity we assume $L(s,\pi)$ holomorphic. We have the relation
$$
\int_1^\infty A_0(x) x^{-s}\frac{dx}{x} = \frac{L(s)}{s}
$$
$(Re(s)>\nu)$. Suppose $s=\sigma+it$, $\sigma>\nu$. This can be written as
$$
\int_0^\infty A_0(e^X) e^{-\sigma X} e^{-itX}dX.
$$
The function $A_0(e^X) e^{-\sigma  X}$ is exponentially decreasing, and therefore this Fourier transform is a $C^0$ function of $t$; this implies that $\mu(\sigma)\le 1$.

The convexity estimate for $\mu(\sigma)$ is $\frac{m}{2}(1-\sigma)$. Thus $\mu(\sigma)\le 1$ beats the convexity estimate if $\sigma<1-\frac{2}{m}$.

Friedlander and Iwaniec have shown that the exponent $\nu=1-\frac{2}{m}$ could be obtained for $L(s) = L(s,\pi_1)L(s,\pi_2)$ and $m=m_1+m_2, m_1,m_2\ge 2$.
See \cite[\S 3]{FI}\footnote{They do not, however, give the proof. It has been explained to the author by Friedlander and Brumley.}. We see that any improvement on $\nu=1-\frac{2}{m}$ would lead to subconvexity (in the $t$--aspect) quite generally for $GL(m)$ (and not only over $\Q$.) \footnote{For $GL(n)$ over $\Q$ this has been recently proved by Nelson \cite{Nel} } If we assume the optimum value conjectured by them, $\nu=\frac{1}{2}-\frac{1}{2m}$, we obtain the following approximation of the Lindel\"of Conjecture:
\begin{equation} 
\mu
(\sigma) \le \frac{m}{2}-\frac{m(m-2)}{m-1}\sigma \qquad (\sigma\le \frac{1}{2}-\frac{1}{2m})
\end{equation}
\begin{equation}
\mu
(\sigma) \le \frac{2m}{m+1}(1-\sigma) \qquad (\sigma\ge \frac{1}{2}-\frac{1}{2m}).
\end{equation}

The restriction to two factors is not one when subconvexity is concerned, because the estimate for $L(s,\pi\times \pi)$ implies one for $L(s,\pi)$.




\vskip6mm

Laurent Clozel

Math\'{e}matiques, B\^{a}timent 307

Universit\'{e} Paris-Saclay

91405 Orsay Cedex

France

\end{document}